\magnification=\magstep1
\dimen100=\hsize
\parskip= 6pt

\font\ninerm=cmr10 at 10truept
\font\eightrm=cmr8
\font\sc=cmcsc10

\font\tenmsy=msbm10
\font\sevenmsy=msbm7
\font\fivemsy=msbm5
\newfam\msyfam
\def\msy{\fam\msyfam\tenmsy}
\textfont\msyfam=\tenmsy
\scriptfont\msyfam=\sevenmsy
\scriptscriptfont\msyfam=\fivemsy

\def\bbc{{\msy C}}

\def\za{\vrule height6pt width4pt depth1pt}

\font\aa=eufm10

\def\Got#1{\hbox{\aa#1}}

\def\Ra*{(R_{a})_{*}}
\def\ada-1{ad_{a^{-1}}~}
\def\p-1u{\pi^{-1}(U)}
\def\gg{{\Got g}}
\def\gp{{\Got p}}

\def\gso{{\Got so}}

\def\gh{{\Got h}}

\def\gd{{\Got d}}

\def\liet{{ \Got t}}
\def \mg {{\cal M}_G}
\def \tmg {T^*{\cal M}_G}
\def \tmgzero {T^*{\cal M}_{G_0}}
\def \mgl{{\cal M}_{GL(N)}}
\def \tmgl {T^*{\cal M}_{GL(N)}}

\def \mred{(\tmg)_{\rm red}}

\font\svtnrm=cmr17

\centerline{\svtnrm Calogero-Moser systems and Hitchin systems} 
\bigskip
\centerline{\sc  J.C. Hurtubise and E. Markman}
\footnote{}{\ninerm The first author of this article would like to thank 
NSERC and FCAR for their support}
\footnote{}{\ninerm The second author was partially supported by NSF grant
DMS-9802532}
\bigskip
\centerline{\vbox{\hsize = 5.85truein
\baselineskip = 12.5truept
\eightrm
{\sc Abstract}: 
We exhibit the elliptic 
Calogero-Moser system as a Hitchin system of $G$-principal 
Higgs pairs. The group $G$, though naturally associated to any root system, 
is not semi-simple. We then interpret the Lax pairs with spectral parameter 
of [dP1]  and [BSC1] in terms of equivariant embeddings 
of the Hitchin system of $G$ 
into that of $GL(N)$. 
}}

{\bf 1. Introduction} 

The Calogero-Moser Hamiltonian system must be one of the most thoroughly 
studied 
Hamiltonian systems, yet many aspects of its geometry remain quite mysterious.

One can associate Calogero-Moser systems to any root system ${\cal R}$ on the 
Lie algebra $\gh$ of a torus $H$ of dimension $r$ and to any elliptic 
curve $\Sigma$ (as well as to the limiting cases of rational nodal or cusp 
curves, the 
``rational" and ``trigonometric'' cases  respectively). In canonical coordinates
$(x,p) \in \gh\times \gh^*$ the system is very simple, and is given by the 
Hamiltonian
$$CM = p\cdot p +  \sum_{\alpha\in {\cal R}} m_{|\alpha|}\gp(\alpha(x)).
\eqno (1.1)$$
where $\gp(x) $ is the Weierstrass $\gp$-function, and the $m_{|\alpha|}$
 are constants depending only on the norm of the root $\alpha$. In the rational case, one replaces $\gp$ 
by the function $x^{-2}$, and in the trigonometric case, by the function $sin(x)^{-2}$.

The system was of course obtained in a step by step fashion from the rational and 
trigonometric 
$Sl(n)$ case, by various people (for a survey, see [OP]), 
who in particular  noticed that one could
 replace the linear functions 
$(x_i-x_j)$ occuring in the original $SL(N)$ case 
by the roots of more general root systems, while maintaining integrability.
 The presence 
of root systems naturally suggests that the  Calogero-Moser systems have some geometric origin,
tied to Lie groups.

In particular, when one is discussing ties between Lie groups or algebras and integrable 
systems, one is immediately led to look for  Lax pairs $\dot L= [M,L]$, and indeed
this is the way much of the work on the Calogero Moser system has progressed,
e.g. with [OP, K] and more recently [dP1, dP2, dP3, BCS1, BCS2], 
so that one now has
 Lax pairs
with a spectral parameter for all of the Calogero-Moser systems. Nevertheless,
 there are several mysterious aspects to many of these Lax
pairs, and in general, there is a lack of concordance between the Lax pairs
and the geometry of the group: 

\item{-} The first is that, for the most part, the Lax
matrices $L$ are not in the Lie algebra
 of the root systems, though they occasionally occur in a symmetric space 
construction associated to the Lie algebra [OP].
Often, however,  they are   in some $Gl(N)$, where $N$  is not even a dimension
of a non-trivial representation of the group.

\item{-} As was pointed out in [Do], while there is a manifest invariance of 
the Calogero-Moser Hamiltonian under the action of the Weyl group, there is 
on the other hand, for
 most semisimple Lie groups,  no Weyl invariant coadjoint orbits.
Such orbits   would   usually be an essential ingredient for a suitable geometric version
of the Calogero-Moser systems. Also, in general, there are no orbits of dimension
twice the dimension of the torus, which one would also want.

\item{-} Finally, one has a Calogero-Moser system for the root systems $BC_n$,
and these do not even correspond to groups.

There is one case for which a satisfactory geometric version of the Calogero-Moser
 system
exists, that of $SL(N, \bbc)$ (see, e.g. [Do]). In this case, one finds that the Calogero-Moser system 
is a generalized Hitchin system over a moduli space of stable pairs over the elliptic 
curve. That is, one considers the moduli space of pairs $(E,\phi)$, where
\item{-} E is a 
 rank $N$  degree 0 bundle  on $\Sigma$ with trivial determinant, and
\item {-}  $\phi$ is a section of $End(E)\otimes K_\Sigma$
with a simple pole at the origin whose residue is a conjugate
of  $m\cdot diag(1,1,1,..1,-N+1)$). 

This realises the phase space in a natural way:
the Hamiltonians are the coefficients of the equation of the
`` spectral curve'' of $\phi$; and one has natural compactifications
of the level sets of the Calogero-Moser hamiltonians as Jacobians of the spectral
curves.

The following note attempts to explain some of the geometry of the 
Calogero-Moser systems for an arbitrary root system ${\cal R}$
in terms of the geometry of a modified Hitchin system. The departure from previous work
 is 
that we do not use as structure group 
the semi-simple group associated to the root system, but rather 
a group which one can construct for any root system,
whose connected component of the identity   is
 the semi-direct product of the
torus with the sum of the root spaces. The Weyl group acts on these,
and  this allows one to construct in a natural
way for any root system
some Weyl invariant coadjoint orbits of the correct dimension, to which one can associate Hitchin
systems, over which the Calogero Moser Hamiltonian  appears naturally. 
  The Lax pairs with 
spectral parameter of [dP1, BCS2] appear in a natural way from  
 embeddings of the Lie algebras of our groups into $Gl(V)$, where $V$ 
is a sum of weight spaces
invariant under the Weyl group $W$ of the root system.
These embeddings are not homomorphisms, but are invariant under the torus and
the Weyl group action.   

While this in some way clears up some of the mystery surrounding 
the Calogero-Moser systems, and in particular addresses the three 
facts outlined above, there are several aspects that are still 
unexplained: the first is that while the Calogero-Moser Hamiltonian
occurs naturally, the other commuting Hamiltonians have no easy
interpretation in our setup and only 
seem to occur naturally after embedding into 
 the Hitchin system for $Gl(V)$. The other 
remaining question is understanding the compactifications of the 
level sets of the Hamiltonians (as Abelian varieties). This would require an
 enlargement of our phase space.
 
Section 2 of this paper is devoted to recalling certain facts about elliptic curves; 
section 3 is similarily devoted to the required facts on 
generalized Hitchin systems.
In section 4, after discussing extensions of the Weyl group by the torus, 
we introduce the group which interests us, and discuss its properties. The next section
is devoted to Hitchin systems with this group as structure group, and we 
show how the Calogero-Moser systems arise. The sixth section is devoted to a discussion 
of how this systems embeds into the Hitchin systems for $Gl(V)$, giving the 
Lax pairs of [dP1], [BCS2].

{\it Acknowledgments:} We would like to thank Ron Donagi and Jim Humphreys
for helpful discussions.

{\bf 2. Line bundles on an elliptic curve} 

Let $\Lambda $ be a non-degenerate lattice in $\bbc$ with generators $2\omega_1,
 2\omega_2$ and let
$\Sigma = \bbc/\Lambda$ be the corresponding elliptic curve. Denote the origin by $p_0$.
We have on $\bbc$ the standard elliptic functions $\sigma(z), \zeta(z)$
with expansions at $z=0$
$$\eqalign { \sigma(z) &= z + O(z^5)\cr \zeta(z) & = {1\over z} + O(z^3),\cr}\eqno(2.1)$$
and periodicity relations
$$\eqalign{\sigma (z+2\omega_i) &= -\sigma(z) exp(2\eta_i(z+2\omega_i)),\cr
\zeta (z+2\omega_i) &= \zeta(z) + 2 \eta_i,}\eqno(2.2)$$
with $\eta_i = \zeta(\omega_i)$. We have 
$${d\over dz}{\rm log}(\sigma(z)) = \zeta(z), 
\quad {d\over dz}\zeta(z) = -\gp(z),$$
where $\gp(z)$ is the standard Weierstrass $\gp$-function. 
 From the periodicity relations, one has that the function 
$$\rho^1(x,z) = {\sigma(z-x)\over \sigma(z)\sigma(x)} e^{x\zeta(z)}\eqno(2.3)$$
is well defined on the elliptic curve with parameter $z$, 
with an essential singularity at the origin, and a single zero at $x=z$.
If we set 
$$\rho^0(x,z) =e^{-x\zeta(z)} \rho^1(x,z),\eqno(2.4)$$
we find that $\rho^0$ has a single pole in $z$ at the origin.
Covering $\Sigma$ by $U_1 = \Sigma-$ (origin), $U_0 =$ disk around origin, we can 
reinterpret the relation $(2.4)$ as saying that $\rho^1, \rho^0$ define 
a section of the line bundle $L_x$ with transition function $e^{-x\zeta(z)}$;
this section has  a single pole at the origin. Let $p_x$ be the point $z=x$ on $\Sigma$; $L_x$  corresponds to the divisor
$p_x - p_0$.

There is another way of representing sections of the line bundle $L_x$, which 
is as functions $f$ on $\bbc$ satisfying  automorphy relations: 
$$f(z+2\omega_i) = f(z)exp(-2\eta_ix).\eqno(2.5)$$
In this way, the function $\rho^0(x,z)$ also  represents a section of $L_x$ with a simple pole at the origin.

For use later on, we note that 
$$ \zeta(z)  \rho^0(x,z)= e^{-x\zeta(z)}{\partial\rho^1(x,z)\over \partial x}
 - {\partial\rho^0(x,z)\over \partial x}.
\eqno(2.6)$$
As a function of $z$, we have  expansions
$$\eqalign{\rho^0(x,z)&= {-1\over z} + {\sigma'(x)\over \sigma(x)} + O(z),\cr
{d\over dx} \rho^0(x,z) &= {d^2 \over dx^2} ({\rm log}(\sigma(x)) +O(z),\cr}\eqno{(2.7)}$$
We have 
$$\rho^0(x,z)\rho^0(-x,z)= \rho^1(x,z)\rho^1(-x,z) = \gp(z) - \gp(x).\eqno (2.8)$$

Finally, consider
$$\eqalign {I_n: \Sigma& \rightarrow \Sigma\cr
z&\mapsto nz,}\eqno(2.9)$$
and let $I_n^*$ denote the induced action on line bundles;
we note that  $I^*_n(L_x) = L_{nx}$. The pull-back $I^*_n(\rho(x,z)
= \rho(x,nz)$ has automorphy factors $exp(-2n\eta_ix)$, and so one can
represent a section by $\rho^0(x,nz)$ in the trivialization given above:
$$\rho^0(x,nz) =e^{-nx\zeta(z)} \hat\rho^1(x,z)\eqno(2.10)$$
for a suitable function $\hat\rho^1$ on $U_1$; one has
$$n\zeta(z)\rho^0(x,nz) = e^{-nx\zeta(z)}{\partial\hat\rho^1\over \partial x}(x,nz)
 - {\partial\rho^0\over \partial x}(x,nz).\eqno(2.11)$$

{\bf 3. Generalised Hitchin systems on an elliptic curve}

Let $ G$ be a complex Lie group. Following [Ma], 
we will consider the moduli space $\mg$  
 of pairs
\medskip

\centerline{ ($G$-bundle $P_G$ on $\Sigma$, trivialisation $tr$ of $P_G$ at
 $p_0$)}

 \noindent and 
its cotangent bundle $\tmg$. For $P_G\in \mg$, let 
$P_\gg$ be the adjoint bundle associated 
to $P_G$, and $P_{\gg^*}$ the associated coadjoint bundle. 
 For any vector bundle $V$, set $V(p_0)= V\otimes {\cal O}(p_0),
V(-p_0)= V\otimes {\cal O}(-p_0)$.
The fibre of
 $\tmg\rightarrow \mg$ at $P_G$  is the 
vector space $H^0(\Sigma, P_{\gg^*}(p_0))$, (the canonical bundle
of $\Sigma$ is trivial) so that the $\tmg$ is a
 space of triples 
\medskip

\centerline{ ($G$-bundle $P_G$ on $\Sigma$, trivialisation $tr$ at
$p_0$, section $\phi$ of $P_{\gg^*}$ with a   pole at $p_0$)}

  Dually, the tangent space to 
$\mg$ is $H^1(\Sigma,P_\gg(-p_0))$. 
 We have, at $(P_G, T)$, the exact sequence 
for $T(\tmg)$:
$$0\rightarrow H^0(\Sigma, P_{\gg^*}(p_0)) \rightarrow T(\tmg) \rightarrow 
H^1(\Sigma,P_\gg(-p_0))\rightarrow 0. \eqno (3.1)$$

The group $G$ acts naturally on the trivialisations, and so acts symplectically
on $\tmg$. The moment map for this action is simply the residue of $\phi$ at 
$p_0$, expressed in the trivialisation $tr$.
One can take symplectic reductions of $\tmg$ under this action, and obtain reduced
moduli spaces $\mred$.

Let $F$ be a homogeneous invariant function on $\gg^* $ of degree $n$, and 
$\omega$ an element of $H^1(\Sigma, K_\Sigma^{-n+1}(-np_0)) \simeq H^1(\Sigma, {\cal O}
 (-np_0))$. Applying $F$ to $\phi$, one obtains an element 
$F(\phi)$ of $H^0(\Sigma, {\cal O}(np_0))$, and so one can define the Hamiltonian $F_\omega$  
on $\tmg$ by
$F_\omega(P_G, \phi) = \{F(\phi),\omega\}$ where $\{,\}$ denotes the Serre duality
pairing. These Hamiltonians descend to the reduced spaces $\mred$, and define 
there an integrable system, the generalized Hitchin system.

{\it Explicit formulae}

We can cover $\Sigma$ by opens $U_0, U_1$ as above, and choose trivialisations
on these opens, with the one on $U_0$ compatible with the trivialisation $tr$.
Let $T = T_{1,0}$ be the corresponding transition function
 $U_0\cap U_1\rightarrow 
G$ for $\Sigma$; sections of $H^0(\Sigma, P_{\gg^*}(p_0))$ can be represented as functions 
$\phi^i:U_i\rightarrow \gg^*$, with $\phi^1$ holomorphic on $U_1$ 
and $\phi^0$ meromorphic on $U_0$ with only one
simple pole at the origin, and $\phi^1 = ad^*(T) \phi^0$ on $U_0\cap U_1$.

We would like to split the sequence (3.1). Represent a one parameter family  of elements
$(P_G(t), tr(t),
\phi(t)) $ of $\tmg$  by $(T(t),\phi^0(t), \phi^1(t))$, with 
$ \phi^1(t)) = Ad^*(T(t)) \phi^0(t)$. At $t= 0$,  the corresponding tangent
vectors are given by $ \dot v = T^{-1} \dot T, \dot \phi^0, \dot \phi^1$,
with 
$$\dot \phi^1 = Ad^*(T)[ (ad^*(\dot v)(\phi^0)) + \dot \phi^0].$$
 Let $$<\ ,\ >:H^0(\Sigma, P_{\gg^*}(p_0))\times H^1(\Sigma,P_\gg(-p_0)\rightarrow \bbc$$
denote the Serre duality pairing; explicitly, it is defined
by $$< \dot \phi,  \dot v> = res_{p_0}(\dot \phi^0\cdot\dot v)$$
At a point $P_G$ of $\mg$, choose a transition function $T= T_{10}: U_0\cap U_1
\rightarrow G$. Let us choose a vector space $V$ of cocycles mapping isomorphically
to $H^1(\Sigma,P_\gg(-p_0)).$ One can split (3.1) by taking for each 
$\dot v$ in $V\simeq H^1(\Sigma,P_\gg(-p_0)$ the vector
 $(\dot v, \dot \phi_{norm})$
such that the pairing of $\dot \phi_{norm}$ with elements of $V$ is zero.
More generally, for any section $a$ of $P_{\gg^*}(p_0)$ over $U_0$, 
let $a^\&$ denote the element of $H^0(\Sigma, P_{\gg^*}(p_0))$
whose pairing with elements of $V$ is the same as that of $a$.
One then has the isomorphism
$$ \eqalign{ T(\tmg) &\rightarrow   H^1(\Sigma,P_\gg(-p_0) )\oplus
H^0(\Sigma, P_{\gg^*}(p_0))\cr
(\dot v, \dot \phi)&\mapsto (\dot v,  (\dot \phi)^\&).\cr} \eqno (3.2)$$

{\sc Proposition (3.3) }\tensl Under this isomorphism, the symplectic form on $\tmg$
becomes, at $(P_G,tr, \phi)$:
$$\Omega( (v', \phi'^\&), (\dot v, \dot\phi^\&)) = 
<v',\dot \phi ^\&> - (\dot v, \phi'^\&> + <[v',\dot v ],\phi>.\eqno (3.4)$$\tenrm

{\sc Proof:} One can parametrise $\mg$ locally by $V$; indeed, 
if $T$ is a transition matrix
for $P$, one has in a neighbourhood of the origin 
a map $V\rightarrow \mg$ obtained by associating to the cocycle $v$ the transition matrix
$T\cdot exp(v)$; this in turn defines a map 
$\rho: V\times V^* = T^*V\rightarrow \tmg$,  which preserves the symplectic
form. With respect to  the splitting (3.2), the 
differential $d \rho$ at the origin is 
$$d\rho (\dot v,\dot \phi) \rightarrow (\dot v, \dot \phi + 
 {1\over 2}(ad^*_{\dot v} \phi)^\&).\eqno(3.5)$$
Substituting in the standard expression for the symplectic form on 
a cotangent bundle  gives (3.4). \hfill\za

The explicit action of the  group $G$   on $T^*(\mg)$ is given by 
$$g(T, \phi^0, \phi^1) = (Ad_g(T), Ad^*_{g^{-1}}(\phi^0), \phi^1),
$$
and the moment map for this action is $res_{p_0}(\phi^0)$.

 From (3.4), we can compute the Hamiltonian vector fields associated to
$F_\omega$:
$$(\dot v, \dot \phi) = ( df\cdot\omega, 0).\eqno (3.6)$$
In other words, the Higgs field part stays as is, but the bundle varies. That
 this is possible is due to the invariance of the function: $ad^*(df)(\phi) = 0$.
We can write an equivalent version of the flow by modifying the transition 
function  by a coboundary:
that is,  we are allowed to modify our trivialisations of the bundle over the  $U_i$,
as long as on $U_0$ the trivialisation is not changed over $p_0$.  Thus, if 
$g_i\in H^0(U_i, P_{\gg^*})$, with $g_0(0) = 0$, we have the equivalent version of the flow:
$$(\dot v, \dot \phi^0, \dot \phi^1) = ( df\cdot\omega+ Ad(T^{-1})(g_1) - g_0, \ 
ad^*(g_0)\phi^0,\ ad^*(g_1)\phi^1).\eqno (3.7)$$
Similarly, on the reduced space, one can modify the flow by a coboundary, but now with 
$g_0(0)$ arbitrary.

{\bf 4. A group associated to a root system}

In this section, we will define the structure group which we will  use
for our Calogero-Moser systems. 
It  will be associated to the root system ${\cal R}$ acting on
the Lie algebra $\gh$ of  torus $H= (\bbc^*)^r$.  

We begin however with a discussion which shows
 that in some sense, passing to a new group is necessary.  
The symplectic reduction leading to  the generalized Hitchin system
$\mred$ depends on a choice of a coadjoint orbit of $G$. In our case of
an elliptic base curve, the 
dimension of $\mred$ is equal to the dimension of the coadjoint orbit. 
We are thus looking for a group, related to the rank $r$ root system, 
and admitting a $2r$-dimensional coadjoint orbit; $2r$ being the dimension of 
the Calogero-Moser system. When the semi-simple group is $SL_n$, 
the coadjoint orbit of ${\rm diag}(1,1, \dots , 1, 1-n)$ is $2n-2$ dimensional.
A general semi-simple group $G_{ss}$ does not have a $2r$-dimensional 
coadjoint orbit. There is however a group $G_0$, naturally associated to 
$G_{ss}$, which does admit $2r$-dimensional coadjoint orbits. 
We consider the group ${\cal O}_0(G_{ss})$ 
of germs at the origin of maps $\bbc\rightarrow G_{ss}$, 
and let $V$ be the subgroup of 
germs of the form $g(z) = h + zg_1 +z^2g_2+...$, \ $h\in H$, 
and $V'$ be the subgroup 
of $V$ of germs  of the form $g(z) = {\rm Id} + zh_1 + z^2g_2+...$, \ 
$h_1\in \gh$. 
Then $G_0= V/V'$. 
$G_0$ is the semi-direct product of $H$ with the direct sum
of the root spaces. The coadjoint orbits of $G_0$ are analyzed in 
Proposition (4.14). 

 There is a natural extension $N(G_0)$ of the Weyl group $W$ by $G_0$. 
Simply consider germs with a leading coefficient in the normalizer $N$ of 
$H$ in $G_{ss}$. $N(G_0)$ acts on $\gg_0^*$ via the coadjoint action. 
We encounter another difficulty: 
For a general semi-simple group, there does not exist any $W$-invariant 
$2r$-dimensional coadjoint orbit in $\gg_0^*$ (i.e., one which is also an 
$N(G_0)$ orbit). Proposition (4.14) shows that 
for such an orbit to exist, we need a non-trivial 
$W$-invariant $H$-orbit in the direct sum
of the root spaces. When the group $G$ is $SL_n$, 
such an orbit is obtained by intersecting the  direct sum
of the root spaces with the coadjoint orbit of 
${\rm diag}(1,1, \dots , 1, 1-n)$. 
More generally, we  relate the existence of $W$-invariant $H$-orbits 
to splittings of the short exact sequence 
$$0\rightarrow H\rightarrow N\rightarrow W\rightarrow 0,\eqno(4.1)$$
 where $N$ is any extension of $W$ by $H$ 
such that the conjugation in $N$ induces on $H$ the standard $W$-module 
structure. 
Let $V$ be an $N$ representation and $R$ a non-zero 
$W$-invariant $H$-orbit in $V$. 
Given a character $\alpha$ of $H$, denote by $N_\alpha$ and $W_\alpha$
its stabilizers in $N$ and $W$.

{\sc Lemma (4.2) }\tensl
\item {1)}
The stabilizer $Stab(\xi)$ of every element $\xi\in R$ intersects $H$
in a fixed normal subgroup $Stab_0\subset H$. 
\item{2)} Let
$$0\rightarrow \overline{H}\rightarrow
 \overline{N}\rightarrow W\rightarrow 0\eqno(4.3)$$
be the quotient of (4.1)
by $Stab_0$. Then the stabilizer $Stab_{\overline{N}}(\xi)$ of every 
$\xi \in R$ projects isomorphically onto $W$. In particular, 
(4.3) splits. 
\item{3)}
If $V$ is an irreducible representation of $N$, and the $W$-invariant 
$H$-orbit $R$ is not the zero orbit, and 
$\iota_\xi : W \hookrightarrow \overline{N}$ is the splitting provided by 
$\xi\in R$, then the representation 
$\iota_\xi^{-1}(V)$ of $W$ is equivalent to 
$Ind_{W_\alpha}^W(1)$ for any character $\alpha$ of $H$ in $V$. 
Consequently, 
we obtain a characterization of $V$ as a representation of $N$:
$V$ is equivalent to the pullback to $N$ of 
$Ind_{\overline{N}_\alpha}^{\overline{N}}{\alpha}$ where ${\alpha}$
is the unique character of $\overline{N}_\alpha$ which restricts to the 
trivial character of $W_\alpha$ and to the character ${\alpha}$ of 
$\overline{H}$.

\tenrm

\noindent
{\sc Proof:}
1) We have the equality $Stab(n\cdot \xi) = n\cdot Stab(\xi)\cdot n^{-1}$. 
Since $R$ is a $H$-orbit, every two stabilizers of elements
in $R$ are conjugate by an element of $H$. 
$Stab(\xi)\cap H$ is a fixed group $Stab_0$ as  $H$ is commutative. 
$Stab_0$ is a normal subgroup of $N$ since
$$[n\cdot Stab_0\cdot n^{-1}\ = \ n\cdot Stab(\xi)\cdot n^{-1}\cap H \ = \
Stab(n\cdot\xi)\cap H \ = \ Stab_0.
 ] $$

2) Let $n_w$ be an element of $N$ mapping to $w\in W$.
Choose an element $\xi \in R$. Since 
$R$ is also an $N$-orbit, 
there exists an element $a\in H$ such that $n_w\cdot \xi = a\cdot \xi$.
Thus, $Stab(\xi)$ maps onto $W$. If follows that 
the stabilizer $Stab_{\overline{N}}(\xi)$ 
in $\overline{N}$ surjects onto $W$. The homomorphism 
$Stab_{\overline{N}}(\xi)\rightarrow W$ is injective because
$Stab_{\overline{N}}(\xi)\cap \overline{H} = (1)$. 

3) Let $\alpha$ be any character of $H$ with positive multiplicity in $V$, 
$V_\alpha$ the corresponding subspace, 
and $\xi$ an element of $R$. 
Since $V$ is irreducible, $\xi_\alpha \neq 0$. 
Since $\iota_\xi(W)$ is the stabilizer of $\xi$, 
the line spanned by $\xi_\alpha$ in $V_\alpha$ is the
trivial character of $W_\alpha$. 
It follows that the direct sum $V'$ of all the translates of 
$span\{\xi_\alpha\}$ 
by $\iota_\xi(W)$ is a sub-representation of $V$. 
The irreducibility of $V$ implies that $V_\alpha=span\{\xi_\alpha\}$ and 
$V'=V$. 
It follows that, as an $N$-module, $V$ is the induced representation 
$Ind_{\overline{N}_\alpha}^{\overline{N}}V_\alpha$. 
The equivalence $V\cong Ind_{W_\alpha}^W(1)$ of $W$-modules 
follows. Note that $V$ need not be irreducible as
a $W$-representation.
\hfill\za

The Lemma specifies two obstructions to the existence of a 
non-trivial $W$-invariant $H$-orbit in $\gg^*$ for a simple Lie algebra.
The first obstruction is the extension class of (4.1).
In particular, considering the (co)-adjoint representation, 
the list of simple groups of adjoint type for which the exact sequence
(4.1) for the normaliser does not split is: 
$SP(n)$, $F_4$, $E_6$, $E_7$, $E_8$ (mod centers). 
 The second obstruction, condition (3) in the Lemma, rules out the 
existence of a $W$-invariant $H$-orbit in $\gg^*$ for Lie algebras of type 
$D_n$ (and in the long root representation of type $B_n$) 
even though (4.1) splits. 

{\sc Example:}   The exact sequence (4.1) splits for $SO(2n)$ 
and $W$ embeds in $N$ as a subgroup of the group of permutation matrices. 
Identify the Lie algebra
$$\gso(2n)=\left\{ 
\left(\matrix{m & n \cr p & q }\right) \ : \ q=-m^t, \ n^t=-n, \ 
{\rm and} \ p^t=-p \right\}, 
$$
consider the Cartan $span\{e_{i,i}-e_{n+i,n+i}\}$ and let $\alpha$ be 
the root $\epsilon_i-\epsilon_j$ corresponding to 
$\gg_\alpha=span\{e_{i,j}-e_{n+j,n+i}\}$, \ $i\neq j$. 
The matrix of the permutation 
$\sigma:=(i,n+j)(j,n+i)$ belongs to the stabilizer 
$W_\alpha$ of the point $\alpha$ in the root lattice. However, 
$\sigma$ acts by multiplication by $-1$ on $\gg_\alpha$. In particular, 
$\gg_\alpha$ is not the trivial character of $W_\alpha$ and 
condition (3) of the Lemma is not satisfied. 

 In order to circumvent the first obstruction, 
instead of the normaliser 
$N(H)$ of the torus, we will consider the semi-direct product $N'$
of the torus and the Weyl group:
$$0\rightarrow H\rightarrow N'\rightarrow W\rightarrow 0.\eqno(4.4)$$

We now define our structure group. 
For any set of  weights $w$ which is invariant under the Weyl group, 
 part (3) of Lemma (4.2) determines 
a representation of $N'$ on the associated sum of weight spaces 
$V=\oplus \bbc_w$: indeed, if we choose a basis
element for each weight space $\bbc_w$, the Weyl group acts simply by 
permuting these basis elements, 
while  the torus acts in the natural way on each weight space.
This holds in particular for the roots $\alpha: \gh\rightarrow \bbc$. 
We define $G$ to be the semi direct product
$$\oplus_{\alpha=1}^n\bbc_\alpha\rightarrow G\rightarrow N'.\eqno(4.5)$$
The connected component of the identity is the group 
$G_0$ discussed above. $G_0$ is the semi-direct product
$$\oplus_{\alpha=1}^n\bbc_\alpha\rightarrow G_0\rightarrow H.\eqno(4.6)$$

Given any element of the Lie algebra $\gg$, we can decompose it into its torus 
and root space components; write this decomposition as 
$ \xi = \xi_\gh + \xi_r$.
Similarly, we can decompose an element $a$ of $\gg^*$ as $a_\gh + a_r$.
The choice of the group $G$ is motivated by the following: 

{\sc Proposition (4.7):} \tensl  The $G-Ad^*$-invariant functions on $\gg^*$ 
only depend on the  root space components, and correspond to the 
$N'$-invariant functions on $\sum_\alpha\bbc_\alpha$. 
The generic coadjoint orbit
is $2r$-dimensional, where $r= dim(N')$, and is of the form 
$$
(N'-{\rm orbit \ in} \  \sum_\alpha\bbc_\alpha) \ \ \times \ \ \gh^*.
$$
 Moreover, $\gg^*$ has a $2r$-dimensional connected 
($W$-invariant) coadjoint orbit.
\tenrm

 Proposition (4.7) follows from Proposition (4.14). 
The rest of this section is dedicated to the proof of these two 
Propositions. 

 Let us fix a basis element for each root space $\bbc_\alpha$,
in a $W$-invariant way.
The components of a vector $C\in
\oplus_{\alpha=1}^n\bbc_\alpha$ are then well defined, and  
naturally indexed  by the roots themselves:
$C= (C_\alpha)$. 
 Let $\tilde\alpha:H\rightarrow \bbc^*$ denote the character 
corresponding to $\alpha$, so that $H$ acts on $\bbc_\alpha$ by 
$(h,v)\mapsto \tilde \alpha(h)\cdot v$. Let $C^t\cdot D$ denote the natural scalar product of 
two vectors in $\oplus_{\alpha=1}^n\bbc_\alpha$, and let $C\circ D$ denote the 
componentwise product: $(C\circ D)_\alpha = C_\alpha D_\alpha$. We denote by $I$
 the permutation matrix acting on the   root spaces which permutes the $\alpha$-th
and the $-\alpha$-th root spaces. Finally, we write the action  of $\gh $ 
on $\oplus_{\alpha=1}^n\bbc_\alpha$ as a matrix: 
$$ \sum_i\tau_iA_{i,\alpha} = \alpha(\tau),$$ so that the action of $\tau$ on
 $C$ can be written as $(\tau^t\cdot A)\circ C$.

As a manifold,  $G_0= (\oplus_{\alpha=1}^n\bbc_\alpha)\times  H$.
The product is given by 
$$(C,h) (C',h') = (C +  (exp( log(h)\cdot   A))\circ C', hh')
.\eqno(4.8)$$
The corresponding Lie Bracket on $(\oplus_{\alpha=1}^n\bbc_\alpha) \oplus \gh$ is 
$$[ (\Gamma, \tau),(\Gamma', \tau')] = ( (\tau^t \cdot  A) \circ \Gamma'  -
(\tau'^t\cdot  A) \circ \Gamma', 0).\eqno(4.9)$$ 
There is a pairing on the Lie algebra, identifying $\gg$ with $\gg^*$:
$$<(\Gamma, \tau),(\Gamma', \tau')> = 
\Gamma^t\cdot I\cdot \Gamma' + 
\tau^t\cdot \tau'.\eqno(4.10)$$
 We will use this pairing to describe the coadjoint action:
one then has 
$$\eqalign {<[(\Delta, \sigma), (\Gamma, \tau)],
(\Gamma', \tau')>&=   ( \sigma^t\cdot A\circ \Gamma - 
\tau^t\cdot A\circ \Delta)^t \cdot I \cdot \Gamma' \cr
&= - \Gamma^t\cdot I\cdot ((\sigma^t \cdot A)\circ\Gamma')
 -  \tau^t \cdot A\cdot(\Delta\circ(I\cdot\Gamma')),\cr}
 \eqno(4.11)$$
remembering that $A_{i,-\alpha} = - A_{i,\alpha}$.
Therefore 
$$ad^*_{(\Delta,\sigma)}(\Gamma', \tau') = 
((-\sigma^t \cdot A)\circ\Gamma', -A \cdot(\Delta \circ (I\cdot\Gamma'))).
\eqno (4.12)$$
Similarily, the coadjoint action of an element of the group can 
be written as
$$Ad^*_{(D, exp(\sigma))}(\Gamma',\tau') = 
(exp (-\sigma^t \cdot A)\circ\Gamma', \tau' -A \cdot(D\circ (I\cdot\Gamma'))).
\eqno (4.13)$$

 From (4.13), one has:

{\sc Proposition (4.14):} \tensl  
The $G_0-Ad^*$-invariant functions on $\gg^*$ 
only depend on the  root space components, and correspond to the 
$H$-invariant functions on $\sum_\alpha\bbc_\alpha$. The generic coadjoint orbit
is $2r$-dimensional, where $r= dim(H)$, and is of the form 
($H$-orbit in $\sum_\alpha\bbc_\alpha$) $\times \gh^*$.\tenrm

The pairing (4.10)  is part of a more general family of invariant inner 
products on  
$\gg$: 

{\sc Lemma} (4.15) \tensl Let $D$  be a diagonal matrix such that 
$D_{\alpha,\alpha}=D_{\alpha',\alpha'}$ if $\alpha,\alpha'$ lie in the same Weyl 
group orbit, and let $\delta$ be a constant. The $N'$ invariant pairings on 
$\gg$ are given by 
$$<(\Gamma, \tau),((\Gamma', \tau')>_{D,d} = 
\Gamma^t\cdot D\cdot  I\cdot \Gamma' + 
\delta\tau^t\cdot \tau'.\eqno(4.16)$$
When there is a single Weyl orbit of roots, there is then a two parameter
family of pairings; when there two orbits, there 
is a three parameter family.
\tenrm

{\sc Proof}: The invariance under $H$ forces the $\alpha$-th root space to be 
paired only with the $-\alpha$-th, and $\gh$ to be paired only with itself.
 Further invariance under the Weyl group reduces one to a single choice 
up to scale for $\gh$ and for each Weyl orbit of roots.

\bigskip

{\bf 5. Hitchin systems for $G$}

We now turn to studying (modified) Hitchin systems for our group $G$, over
an elliptic curve $\Sigma$. We begin with 
$$\mg = \{ G-{\rm bundles \ of \ degree \ zero, \ trivialized \ at } \ p_0\},$$ 
then take the cotangent
bundle of this space. We then reduce by the action of the group $G$, at a 
$W$-invariant element of $\gg^*$; this element must then lie in $\gg^*_r$. We will see that this is essentially 
equivalent to reducing by the action of the subgroup $\oplus_\alpha\bbc_\alpha$.
The Calogero-Moser Hamiltonians are then expressed naturally in terms of the 
scalar product  of (4.10) on the reduced space. 

{\it 5.i. $G$-bundles on $\Sigma$}

We begin by giving an explicit description of the moduli of framed 
$G$-bundles of degree 0 on $\Sigma$. 
We first note that under the projection of 
$G$ to $W$, any $G$-bundle defines a $W$-bundle. 
We will consider only the component  of moduli corresponding to trivial 
$W$-bundles, so that we can represent our bundles as 
$G_0$-bundles. The subspace of $G$-bundles which can be represented as $G_0$
bundles is a quotient of the moduli of $G_0$-bundles:
one must quotient out by the action of $W$, since different $G_0$-bundles can be 
the same as $G$-bundles. 
$${\cal M}_{G } = {\cal M}_{G_0}/W. \eqno (5.1)$$

Let  $\overline{{\cal M}}_{G_0}$ be the moduli of 
$G_0$-bundles (without framing). 
To analyse the moduli $\overline{{\cal M}}_{G_0}$, we use the fact that 
the group $G_0$ maps to $H$, and so one has maps 
$$\Pi : \overline{{\cal M}}_{G_0}  \rightarrow \overline{{\cal M}}_H  = Pic^0(\Sigma)^r = 
\Sigma^r, \eqno (5.2)$$
$$\Pi_W: \overline{{\cal M}}_{G_0}/W \rightarrow \overline{{\cal M}}_H/W = 
\Sigma^r/W.\eqno (5.3$$
By a theorem of Looijenga [Lo], the space  $\Sigma^r/W$ 
is a weighted projective space.  
The fiber of (5.1) at $\chi\in \Sigma^r$ is   $\oplus _\alpha 
 H^1(\Sigma,L_{\tilde \alpha(\chi)})$,  where $L_{\tilde \alpha(\chi)}$ 
is the line bundle associated to ${\tilde \alpha(\chi)}$. 
This can be seen by writing out a cocycle explicitly in the semi-direct 
product.
Each   $H^1(\Sigma,L_{\tilde \alpha(\chi)})$ is isomorphic to 
$\bbc$ if $L_{\tilde \alpha(\chi)}$ is trivial and is $(0)$ otherwise. 
  Consequently, one has an open set 
$\overline{{\cal M}}'_{G_0}\subset \overline{{\cal M}}_{G_0}$ 
isomorphic to the 
open set of $\Sigma^r/W$ corresponding to $H$-bundles for which none of the 
$L_{\tilde \alpha(T_\gh)}$ are trivial.

Putting the framings back in, one has that the moduli of framed 
$H$- bundles is the same as the moduli of unframed $H$-bundles, as the
automorphisms act transitively on framings. Consequently, one has
$$\Pi :  {{\cal M}}_{G_0}  \rightarrow \overline{{\cal M}}_H  = Pic^0(\Sigma)^r = 
\Sigma^r, \eqno (5.4)$$
$$\Pi_W:  {{\cal M}}_{G_0}/W \rightarrow \overline{{\cal M}}_H/W = 
\Sigma^r/W.\eqno (5.5).$$
 This time, the fibre is  $\oplus_\alpha H^1(\Sigma,L_{\tilde 
\alpha(\chi)}(-p_0))$. 
Each $H^1(\Sigma,L_{\tilde \alpha(\chi)}(-p_0))$ is isomorphic to $\bbc$. 

{\sc Definition}:  
\tensl 
We say that a $G$-bundle $P_G$ is {\rm special} if $L_{\tilde \alpha(T_\gh)}$
is trivial for some root $\alpha$.
\tenrm

One has an open set 
$$\mg'\subset \mg$$   of framed non-special $G$-bundles.
We will take a reduction of $\mg'$, which will be the space over which 
the Calogero-Moser systems are defined, and indeed, we shall see that 
the reduction  does not 
extend to  the locus
where one of the  $L_{\tilde \alpha(T_\gh)}$ is trivial. 

Explicitly, covering the elliptic curve by $U_0$ = disk around $p_0$, and
 $U_1$ = $\Sigma-p_0$, 
the torus part $T_\gh$ of the transition function $T$
is a  function $T_\gh:U_0\cap U_1\rightarrow H$.
 We can choose these functions 
to be of the form
$$T_\gh  = exp(x\zeta(z)),\ x\in \gh.\eqno(5.6)$$
 The root
space  part $T_r$ can be represented by a vector $M$ of cocycles 
$M_\alpha$ representing 
elements of
$H^1(\Sigma,L_{\tilde \alpha(T_\gh)}(p_0))$.
 For $L_{\tilde \alpha(T_\gh)}$
non trivial,
these cocycles can be taken to be constant functions on $U_0\cap U_1$; when
$L_{\tilde \alpha(T_\gh)}$ is trivial, the constant functions
 correspond to trivial classes, and one must choose a function with a simple 
pole at $p_0$ as generator, for example $\zeta(z)$.


One has as cotangent space to $\mg'$ at a bundle $P_G$
the set of Higgs fields $\phi$ in $H^0(\Sigma, P_{\gg^*} (p_0))$ 
(fields in the associated coadjoint bundle with 
poles at $p_0$). Splitting 
$\phi$ into a root space component and a torus component,
$$ \phi=    \phi_\gh +    \phi_{r},\eqno(5.7)$$
one has that the components $\phi_\alpha$ of  $\phi_r$  have  poles
at the origin only when the line bundle  is not trivial. 
Explicitly, for a bundle with transition functions $(T_r, T_\gh =exp(x\zeta(z)))$,
 one represents 
 $( \phi_r,\phi_\gh)$ in the $U_0$-trivialisation by $(\phi_r^0,\phi_\gh^0)$, 
and in the $U_1$-trivialisation by $(\phi_r^1, \phi_\gh^1)$, with  
$$(\phi_r^1, \phi_\gh^1)= 
(exp(- x^t\cdot A \zeta(z))\circ \phi_r^0
 , \phi_\gh^0 -A\cdot(T_r\circ(I\cdot \phi_r^0))).$$
 Here    $\phi_r^0$ has simple poles at the 
origin; its components  $\phi^0_{ \alpha}$ are simply multiples of the functions
$\rho^0(\alpha(x), z)$ of (2.4). Similarly, the components $\phi_\alpha^1 $ of 
$\phi_r^1 $ are   multiples of the functions
$\rho^1(\alpha(x),z)$.

{\it 5.ii. Reduction}

There is, as in section 3,  an action of $G$ on $\tmg$
by changing the trivialisation at $p_0$.
The moment map for this action is, as we saw in section 3, simply the
residue of the Higgs field $\phi$ at $p_0$, expressed in the trivialization.
We want to reduce  at an element 
$C$ of $\gg^*$ which is $W$ invariant. 
Being $W$-invariant, $C$ must lie in $\gg^*_r$. $W$-invariance 
implies that reduction of 
$\tmg$ by $G$ at $C$ is equivalent to the reduction of $\tmgzero$ by $G_0$
at $C$, then quotienting by $W$.  
The set of  root vectors $\alpha$ splits up into $W$-orbits according to 
their lengths $|\alpha|$, and we choose constants $c_{|\alpha|}$ for each 
length, and set $C = (c_{|\alpha_1|},c_{|\alpha_2|},..,c_{|\alpha_n|})$.
 Note that the coadjoint orbit $G_0\cdot C$ is also a $G$-orbit,
the connected orbit of Proposition (4.7).  

If the constants $c_{|\alpha_i|}$ are non-zero, 
the coadjoint orbit of  the element $C$ is of the form ($H$-orbit in 
$\gg^*_r)\times \gh^*$. Also, the stabiliser of the element $C$ lies 
in the root space part of the group.  
 From these facts and from the expressions for 
the actions, it follows that taking the symplectic quotient by 
 $G_0$ or the quotient by its subgroup $\oplus_\alpha\bbc_\alpha$ gives 
exactly the same result. 

The action of an element $V\in \oplus_\alpha\bbc_\alpha$ on 
$(T,\phi^0)\in\tmg$ is given explicitly by 
$$(T,\phi^0) = ((T_r, T_\gh), (\phi_r^0, \phi_\gh^0 ))\mapsto 
((T_r+V, T_\gh),
 (\phi_r^0, \phi_\gh^0 -A\cdot(V\circ(I\cdot \phi_r^0)).\eqno(5.8)$$
The moment map for this action  
is  res$(\phi_r^0)$.

  Reducing at $C$, we fix to $C$ the residues of $\phi_r^0$, and quotient  by 
the group. Referring to the explicit form of the action, this means that
we can normalise to 
$ T_r = 0$, thus reducing the bundle to the torus, over the locus $\mg'$. 
Our reduced space $\mred$ can thus be thought of as a subspace of the 
unreduced one. This subspace $\mred$ of $\mg$ is characterised by:
$$T_r = 0,\quad {\rm res}(\phi_r^0) = 
{\rm a}\ W-{\rm invariant\ constant}\ 
C.\eqno(5.9)$$

Note that we still have the torus part of the framing  in this 
description of our reduced space.

Explicitly, the elements of $\mred$ are then $H$-bundles with 
transition functions 
$$T = (T_r,T_\gh) = (0, exp(x\zeta(z)),$$
along with Higgs fields whose root space components are, in the $U_0$-trivialisation
$$\phi^0_\alpha = c_{|\alpha|}\rho^0(\alpha(x), z),$$
and whose torus components are constants
$$\phi_\gh^0=\phi_\gh^1= p.$$ 
By (3.3) the functions $x\in\gh,p\in \gh^*$ provide canonical coordinates on 
$\mred$; one must remember that we are restricting to the set 
$\alpha(x)\neq 0, \alpha\in {\cal R}$,
and quotienting out by the  action of  the affine Weyl group 
on $\gh\times \gh^*$.

 Remarks: 
1) In the above discussion we have excluded a Higgs
pair $(P_G,\phi)$ if 
the bundle $P_G$ is special (see section 5.i). 
This exclusion is {\it automatically imposed} by the reduction along the 
coadjoint orbit $G\cdot C$.
The coadjoint bundle $P_{\gg^*}$ fits in the short exact sequence
$$
0\rightarrow 
\gh^* 
\ \rightarrow \ 
P_{\gg^*} 
\ \rightarrow \ 
(\oplus_{\tilde{\alpha}} L_{\tilde{\alpha}})^* 
\ \rightarrow 0.
$$
The moment map sends $(P_G,\phi)$ into $G\cdot C$ if and only if
the residue of $\phi$ projects into the $H$-orbit of $C$ in the fiber
of $(\oplus_{\tilde{\alpha}} L_{\tilde{\alpha}})^*(p_0)$ at $p_0$. 
The triviality of one of the $L_{\tilde{\alpha}}$ rules out the existence 
of such sections in $H^0(\Sigma,P_{\gg^*}(p_0))$. 

2) All the non-special principal $G$-bundles considered have a canonical 
reduction to a principal $N'$-bundle (equivalently, 
a $W$-orbit of reductions to an $H$-bundle). 
This follows from the fact that the global sections of $P_{\gg}$ 
generate a sub-bundle of commutative sub-algebras isomorphic to $\gh$. 

\medskip
{\it 5.iii The Calogero-Moser Hamiltonian}

Let $\omega$ be the  class in $H^1 (\Sigma, {\cal O})$, with 
representative with respect to the cover $U_0,U_1$
$$\omega = \zeta (z).\eqno (5.10)$$ Our Hamiltonian on $\mred$ will  then be the 
W-invariant function on the reduced space: 
$$  CM= res_{p_0}(\omega <\phi,\phi>)\} .\eqno(5.11)$$

Note that the  bilinear form of Lemma 4.15 
gives rise to a canonical bilinear form on $P_{\gg}$ because all 
the non-special principal $G$-bundles we consider have a canonical reduction 
to $N'$. 
The bilinear form sends a Higgs pair $(P_G,\phi)$ to the element 
$<\phi,\phi>$ of $H^0(\Sigma, K_{\Sigma}^{\otimes 2}(2p_0))$. 
$H^0(\Sigma, K_{\Sigma}^{\otimes 2}(2p_0))$ is two
dimensional. If the pair $(P_G,\phi)$ belongs to $\mred$, then
the residue of $\phi$ belongs to our coadjoint orbit $G\cdot C$. Hence, 
 the quadratic residue of $<\phi,\phi>$ is fixed and $<\phi,\phi>$ 
lies in a marked affine line $\ell$ in 
$H^0(\Sigma, K_{\Sigma}^{\otimes 2}(2p_0))$. 
We see that the Calogero-Moser Hamiltonian is determined canonically up to
a choice of an affine linear isomorphism $\ell\cong \bbc$.
The restriction of the 
linear functional (5.11) on $H^0(\Sigma, K_{\Sigma}^{\otimes 2}(2p_0))$ 
provides such an isomorphism. 

One can split $CM$ into a sum $CM_r + CM_\gh$ of a root space piece
$CM_r = \{\omega,<\phi_r,\phi_r>\}$ $ = \{\omega,
 \sum \phi^0_{\alpha}\phi^0_{-\alpha} \}$ and a torus piece 
$CM_\gh =\{\omega,<\phi_\gh,\phi_\gh>\}.$
The relation (2.8) tells us that 
$$\phi^0_{\alpha}\phi^0_{-\alpha}=c_{|\alpha|}^2(\gp(z) -\gp(\alpha(x))$$
and so   the Hamiltonian is 
$$CM = p\cdot p -  \sum_\alpha c_{|\alpha|}^2\gp(\alpha(x)),\eqno (5.12)$$
which is indeed the Calogero-Moser Hamiltonian, setting 
$m_{|\alpha|} = - c_{|\alpha|}^2$.

 Next we provide an explicit formula for the 
vector field of the Calogero-Moser Hamiltonian. The formula (5.18) 
will be needed in Section 6.
The function $CM_r$ is $ad^*$-invariant, and its differential on the
$H^0(\Sigma, P_{\gg^*}(p_0))$-compo\-nent of the tangent space $\tmg$ is
 then 
$$dCM_r= (\omega \phi)^0_r ,\eqno(5.13)$$
where one thinks of $dCM_r$ as an element of $H^1(\Sigma,P_\gg(-p_0) )$
acting on $H^0(\Sigma, P_{\gg^*}(p_0))$. With respect to the splitting (3.2)
the action of $dCM_r$ on the $H^1(\Sigma,P_\gg(-p_0) )$-component of the tangent space
is trivial.
By the considerations of section 2, 
the Hamiltonian vector field in $\tmg$ of $CM_r$ at $(T, \phi^0,\phi^1)$,
 where $T$ is the 
 (torus) transition matrix and $\phi^1= Ad^*(T)(\phi^0)$ 
is the Higgs field, is given by:
$$T^{-1}\dot T = \omega  \phi_r^0, \ \dot \phi^i = 0\eqno(5.14)$$
This takes us out of the normalised form for the reduced space,
 since the transition function no longer lies in the torus. 
Remembering that we are on an elliptic curve,  
the root space components   $\omega \phi_\alpha$ in $H^1(\Sigma, L_\alpha)$
(in the reduced space, we have forgotten the 
root space component of the framings, and so we are in  $H^1(\Sigma, L_\alpha)$
instead of $H^1(\Sigma, L_\alpha) (-p_0))$) can be 
written as coboundaries:
$$(\omega \phi)_{r}= (\omega \phi)^0_{r}+   Ad_ {T^{-1}}(\omega \phi) ^1_{r} 
\eqno(5.15)$$
and  the flow of $CM_r$ 
can be written
$$\eqalign{T^{-1}\dot T&= 0,\cr 
\dot \phi^0 &= (ad^*_{  (\omega \phi)^0_{r} }
\phi^0 )^\& .\cr
}\eqno(5.16)$$
The vector field on $\mred$ corresponding to $CM_\gh= p\cdot p$, with 
respect to the splitting (3.2) is simply 
$$\eqalign{T^{-1}\dot T&=\omega \phi^0_\gh = \omega p,\cr 
\dot \phi^0 &= 0\cr
}\eqno(5.17)$$
and so combining (5.16) and (5.17), one has for the flow of $CM$
$$\eqalign{T^{-1}\dot T&=\omega\phi^0_\gh,\cr 
\dot \phi^0 &=(ad^*_{( \omega \phi ^0_{r})}
\phi^0 )^\&. \cr
} \eqno(5.18)$$

Let us check that the   vector field (5.18 )  is indeed the 
vector field that one obtains from the explicit parametrisation.
 From (2.6), using $\omega= \zeta(z)$ we find that the $\alpha$-th components of the coboundary
decomposition (5.15) satisfy
$$\eqalign{ (\omega \phi)^0_{\alpha}A_{i,\alpha}=& {d \over dx_i} (\phi^0)_\alpha = 
c_{|\alpha|}{d\over dx_i}\rho^0(\alpha(x),z)\cr
},\eqno(5.19)$$
where we decompose $x\in \gh$ into components $x_1,..., x_r$, and $\alpha_i= 
{d\alpha\over dx_i}$ is the corresponding
component of the root $\alpha$. Referring to the formulae (4.12) for the 
coadjoint action and to (2.8), we 
have for the flows:
$$ \eqalign {\dot p_i &= \sum_\alpha   (c_{|\alpha|}{d\over dx_i}\rho^0(\alpha(x),z))
\cdot c_{|\alpha|}\rho^0(-\alpha(x),z)\cr
&= {1\over 2} \sum_\alpha (c_{|\alpha|}^2 {d\over dx_i} \gp(\alpha(x)),\cr
\dot x_i &= p_i.}\eqno(5.20)$$

{\bf 6. Embeddings in $Gl(N,\bbc)$.}

We now give two embeddings of our system into the Hitchin systems for
$Gl(N,\bbc)$   over $\Sigma$,
one corresponding to the Lax pairs of [dP1], the other to those of [BCS2]. 

Let $V= \bbc^N$ be a sum of (integral) weight spaces $\bbc_{\omega_i}, i =1,...
,N$ for the torus $H$ , such that 
the set of roots is Weyl invariant. The weights $w_i$ are maps of $\gh$ to $\bbc$; denote the corresponding 
homomorphism $H\rightarrow \bbc^*$ by $\tilde w_i$.
As for the root spaces, each of these weight spaces should be thought of as having a preferred basis, 
and the bases are
invariant under the Weyl group.

One  has an embedding  $\Xi$ of  the torus $H$ into the diagonal subgroup $D$ 
of $gl(N,\bbc)$; it is given by 
$$\Xi(h)={\rm diag}(\tilde w_1(h),.....\tilde w_N(h)).\eqno(6.1)$$
Let $\xi$ denote the corresponding Lie algebra homomorphism. 

The homomorphism $\Xi$ induces a map $\hat \Xi $ from the space of 
$H$-bundles of degree 0 over $\Sigma$
to the space of $D$-bundles of degree 0 over $\Sigma$, 
where $D$ is the diagonal subgroup of $Gl(N)$. 
The space of $H$-bundles, as we saw, can be parametrised with some redundancy
by $\gh$; for $h\in \gh$, the corresponding $D$-bundle $\hat \Xi(E)$ is a sum 
of line bundles $\oplus_i L_{w_i(h)}$. 
The bundle $End(\hat\Xi(E))$ is then a sum of 
line bundles $\oplus_{i,j}L_{w_i(h)-w_j(h)}$. 
The differences  $w_i-w_j$ are sums of roots, for $w_i, w_j$ in the 
same orbit. We note that the space of  $Gl(N,\bbc)$-bundles 
is essentially the  finite quotient of
the space of $D$-bundles  by the Weyl group of $Gl(N)$,  
as the generic $Gl(N,\bbc)$-bundle reduces to the torus.

{\bf The d'Hoker-Phong embedding}

The map $\hat \Xi $ extends to a map 
$$\Xi_{\cal M}:\mred\rightarrow \tmgl\eqno (6.2)$$
 of the reduced moduli space $\mred$ into 
the cotangent bundle $\tmgl$ of the space of $Gl(N)$- bundles with level 
structure at the point $p_0$ of $\Sigma$. 

  Conceptually, the map $\Xi_{\cal M}$ is determined by an 
$N'$-equivariant  extension of $\xi^*$ to a linear map from $\gg^*$ to $gl(N,\bbc)$. 
Recall that a non-special $G$-bundle admits 
a canonical reduction to a $N'$-bundle. Thus, $\mred$ is also 
a moduli space of pairs 
$(P,\phi)$ where $P$ is a principal $N'$-bundle, $ P_{\gg^*}$ is the vector 
bundle associated via the map of  $\gg^*$ into $gl(N,\bbc)$, 
and $\varphi$ is a section of $P_{\gg^*}\otimes K(p_0)$. The homomorphism 
$\Xi$ extends to a homomorphism $\Xi : N' \rightarrow N(D)$
realizing $gl(N,\bbc)$ as an $N'$ representation. Thus, an $N'$-equivariant
linear map from $\gg^*$ to $gl(N,\bbc)$ gives rise to a map (6.2).

More explicitly, $\mred$ is a space of pairs 
\smallskip

\centerline{($H$-bundle $E$
with $H$-level structure at $p_0$,
section $\phi_G$ of $E( \liet\oplus(\oplus_\alpha\bbc_\alpha))\otimes 
K_\Sigma(p_0)$).}

To this, $\Xi_{\cal M}$ will  associate an element of $\tmgl$.
Such an element is a pair
\smallskip

\centerline{ (rk $N$   bundle $  E_{Gl(N)}$ with level structure at $p_0$,
 section $\phi_{Gl(N)}$ of $End(E_{Gl(N)}) \otimes K_\Sigma(p_0))$.}

The    $Gl(N,\bbc)$-bundle   $  E_{Gl(N)}$  associated by $\Xi_{\cal M}$
to $(E,\phi)$ is simply $\hat\Xi(E)$. We then define the corresponding 
$\phi_{Gl(N)}$. We choose for each pair $(w,w')$ of weights a constant 
$C_{w,w'}$ in   a way that it is invariant under the Weyl group and so that 
$$C_{w,w'}=C_{w',w}.$$
 We then
define a ``shift'' operator for each root $\alpha$
$$(Sh_\alpha)_{w,w'}= \delta_{w-w',\alpha} C_{w,w'},\eqno(6.3)$$
where we index the entries of the matrix by the weights themselves. 
The coefficient $\delta_{w-w',\alpha}$ is the Kronecker $\delta$. 
We then set 
$$\phi_{Gl(N)} = \xi( (\phi_G)_\gh )+ \sum _{\alpha\in {\cal R}}
 (\phi_G)_\alpha Sh_\alpha.\eqno(6.4)$$

 Let 
${\cal CM}$ 
denote the image  $\Xi(\mred)$.
The   space $\mgl$ has dimension $N^2$. Indeed, the space of bundles 
is of dimension $N$:
the generic  $Gl(N)$-bundle on $\Sigma$ reduces to the subgroup $D$ of
 diagonal matrices.
The bundles  have, generically, the group $D$ as automorphisms. 
When one adds in the 
level structure, one adds in $N^2$ parameters, on which the automorphisms act,
reducing one to $N^2-N$ parameters, giving $N+N^2-N= N^2$ parameters in all.
When one considers the Higgs fields $\phi_{Gl(N)}$ in $H^0(\Sigma, 
End(E_{Gl(N)})) \otimes K_\Sigma(p_0))$,
 one has similarily $N^2$ parameters,
giving $2N^2$ parameters for $\tmgl$.
The Calogero-Moser locus  ${\cal CM}$ lies inside $\tmgl$, 
and is of dimension 
$2r$. It is characterised by the fact that the framing is compatible with the reduction
to the diagonal subgroup $\Xi(H) \subset D$ (so that transition functions respecting the trivialisation 
can be chosen diagonal), and the polar parts of the Higgs field $\phi_{Gl(N)}$
are fixed, while its diagonal parts lie in $\xi(\gh)$

More explicitly, for any matrix $A$, let 
$$A = A_d + A_{od}\eqno (6.5)$$ denote the splitting of $A$ into diagonal and off-diagonal
matrices.
One can choose  the transition matrices $M$ for a $GL(N)$-bundles, 
at least generically, to be diagonal, so that
$$ M_d = diag(exp (y_i \zeta(z))),\ M_{od} =  0, \eqno (6.6)$$
where $y_i$ are constant on $\Sigma$. In turn, one represents the
 Higgs fields, which decompose into a sum of sections of line bundles,
by 
$$\eqalign{(\phi_ {Gl(N)} )_d&= diag (q_i),\cr
 ((\phi^0_ {Gl(N)} )_{od})_{w,w'}&= C_{w,w'}(\sum_\alpha \delta_{w-w',\alpha} c_{|\alpha|}) 
\rho^0( \alpha(x),  z), \cr
  }\eqno(6.7)$$
where $q_i$ are constant functions, and $\rho^0$ are the functions of
 (2.4). Note that the $ K_{w,w'}$ are the residues of the section at the origin.
The Calogero-Moser locus  ${\cal CM}$ is given by  constraints
$$ \eqalign{ (1)&\ diag(y_i) \subset \gh,\cr (2)&\ M_{od}=0,\cr 
(3)&\  diag (q_i) \subset \gh,\cr (4)&\ res_0(\phi^0_{Gl(N)})_{w,w'}  = 
C_{w,w'}\sum_\alpha \delta_{w-w',\alpha} c_{|\alpha|}.  } \eqno(6.8)$$

Referring to the explicit form of the symplectic form on
$\tmgl$ (3.3), it follows:

{\sc Proposition  (6.9)}: \tensl  
The embedding $\Xi_{\cal M}$ is symplectic,
and the 
Calogero-Moser Hamiltonian corresponds under the embedding to a multiple of 
$res(\omega tr (\phi^2_{Gl(N)}) )$.\tenrm

We need to compare the flow in the reduced space of $G$-bundles and the 
flow in the cotangent space of $Gl(N)$-bundles with level structure;
 in our $Gl(N)$ moduli, the  Hamiltonian is given, up to a constant,
 by pairing $tr( \phi_{Gl(N)}^2) $ with the cocycle $\omega$ of (5.7); 
the corresponding flow
of $(M, \phi_{Gl(N)})$ is, by (3.6)
$$ M^{-1}\dot M  = \omega  \phi^0_{Gl(N)},\quad  \dot \phi^0_{Gl(N)} = 0.
\eqno (6.10)
$$
Again we split $\omega \phi^0_{Gl(N)}$, 
first into its diagonal and off-diagonal components,
and then write the off diagonal term as a coboundary (function on $U_0$
vanishing at the origin, minus $M^{-1}$(function on $U_1$)  plus
 a constant cocycle:
  $$\omega \phi^0_{Gl(N)} = (\omega \phi^0_{Gl(N)})_d
+ (\omega \phi^0_{Gl(N)})^0_{od} - M^{-1}(\omega \phi^0_{Gl(N)})^1_{od}M 
+(\omega \phi^0_{Gl(N)})^{\rm cst}_{od}\eqno (6.11)$$
referring to (3.7) this transforms the flow (6.12) into the equivalent one:

$$\eqalign { M^{-1}\dot M &= (\omega \phi^0_{Gl(N)})_d +
 (\omega \phi^0_{Gl(N)})^{\rm cst}_{od}, \cr
 (\dot \phi^0_{Gl(N)})  & = [(\omega \phi^0_{Gl(N)})^0_{od},(\phi^0_{Gl(N)}) ],\cr
}.\eqno(6.12)$$
This is not necessarily tangent to the embedded Calogero-Moser 
locus  ${\cal CM}$: 
it does not satisfy the constraints (1), (2) and (3) of (6.7),
but does satisfy the constraint (4).
When one has a symplectic subvariety $V$ of a larger subvariety $W$, one can split
the tangent space of $W$ along $V$ into $TV\oplus (TV)^\perp$, using the symplectic form.
 The Hamiltonian 
vector field of
$H$ along $V$ is simply the projection of the corresponding field in $W$,
 with respect to this splitting. Referring to the formula (3.4), splitting
the diagonal matrices as $\gd= \gh \oplus\gh^\perp$, and letting 
$\pi_\gh:gl(N)\rightarrow \gh,\pi_{\gh^\perp}:gl(N)\rightarrow \gh^\perp$ 
be the ensuing projections, 
the bundle $(T{\cal CM})^\perp$ is given by
$$\eqalign{&\pi_\gh(M^{-1}\dot M) =0\cr
&\pi_\gh(\dot  \phi^0_{Gl(N)} +[M^{-1}\dot M, \phi^0_{Gl(N)}])=0}.\eqno(6.13) $$

 One can make the $\mgl$-flow tangent to the Calogero-Moser locus 
$V$ by adding to it the following vector field, which lies in 
 $T{\cal CM}^\perp$
$$\eqalign{ M^{-1}\dot M &= -(\omega \phi^0_{Gl(N)})^{\rm cst}_{od}\cr
(\dot \phi^0_{Gl(N)})  &= 
a +[\omega \phi^0_{Gl(N)})^{\rm cst}_{od}, (\phi^0_{Gl(N)})]_d},\eqno (6.14)$$
where $a = a(x) $ is a suitable constant (in $z$) in $\gh^\perp$, giving the Calogero-Moser flow:
$$\eqalign{  M^{-1}\dot M &= (\omega\phi^0_{Gl(N)})_d\cr
(\dot \phi^0_{Gl(N)})  &= a(x) +[(\omega \phi^0_{Gl(N)})^0_{od},(\phi^0_{Gl(N)})]
+ 
[\omega \phi^0_{Gl(N)})^{\rm cst}_{od}, (\phi^0_{Gl(N)})]_{d} ,\cr
&= a(x) +[(\omega \phi^0_{Gl(N)})^0_{od}
+ (\omega \phi^0_{Gl(N)})^{\rm cst}_{od},(\phi^0_{Gl(N)}) ]- 
[(\omega \phi^0_{Gl(N)})^{\rm cst}_{od}, (\phi^0_{Gl(N)}) ]_{od} ,\cr}
\eqno (6.15)$$
Indeed, this satisfies the constraints (1) (2) (4) of (6.8), the third constraint 
being given by an appropriate choice of $a\in \gh^\perp\subset \gd$:
$$a(x) = -\pi_{\gh^\perp}([(\omega \phi^0_{Gl(N)})^0_{od}+ 
(\omega \phi^0_{Gl(N)})^{\rm cst}_{od}, (\phi^0_{Gl(N)})]).\eqno (6.16)$$

Let $D'$ be the group of automorphisms of the bundle given by $M$; with respect 
to our trivialisations the automorphisms are represented by constant matrices.
$D'$ includes the group $D$ of diagonal matrices. The action of $D'$ on $T^*\mgl$ is
represented by the vector field: 
$$M^{-1}\dot M = 0, \dot \phi^0_{Gl(N)} = [d',\phi^0_{Gl(N)}],\eqno (6.17)$$
for $d'\in Lie (D')$. 
If $d'$ is diagonal, this vector field lies in ${\cal CM}^\perp$.

We would like to use this vector field to rewrite the flows (6.15) as
$$\eqalign{  M^{-1}\dot M &= (\omega\phi^0_{Gl(N)})_d,\cr
(\dot \phi^0_{Gl(N)}) &= [(\omega \phi^0_{Gl(N)})^0_{od}
+ \omega \phi^0_{Gl(N)})^{\rm cst}_{od} + d'(x),(\phi^0_{Gl(N)})],\cr}
\eqno (6.18)$$
giving a Lax pair with spectral parameter for the flow. This gives the 
constraint
$$\eqalign{[d'(x),(\phi^0_{Gl(N)})] &= 
a- [(\omega \phi^0_{Gl(N)})^{\rm cst}_{od}, 
(\phi^0_{Gl(N)}) ]_{od}.\cr 
}\eqno (6.19)$$
We have, referring to (2.6)-(2.8), 
$$\eqalign{(\phi^0)_{w, w+\alpha}(x,z)&= C_{w, w+\alpha}c_\alpha
(-z^{-1}+  \zeta(\alpha(x)) + O(z)),\cr
 &{\buildrel {\rm def} \over =}R_{w, w+ \alpha} z^{-1} + Q_{w, w+\alpha}(x) +
 O(z).\cr}\eqno(6.20)$$
$$\eqalign{((\omega \phi^0_{Gl(N)})^0_{od} +
+ \omega \phi^0_{Gl(N)})^{\rm cst}_{od})_{w, w+ \alpha}&= C_{w, w+\alpha}c_\alpha
{d\over d\alpha(x)}(\rho^0(\alpha(x), z)\cr
&= R_{w, w+ \alpha}\gp(\alpha(x)) + O(z)\cr
 &{\buildrel {\rm def} \over =} P_{w, w+\alpha}(x) + O(z).}\eqno(6.21)$$

Recall that we are all along dealing with flows of bundles and of sections 
$\phi_{Gl(N)}$, and in particular, that sections are determined by their leading order terms 
at $z= 0$. This gives necessary and sufficient  algebraic constraints for $d'= d'(x)$:
$$[P(x), R]_{od} = [d'(x), R]_{od},\eqno(6.22)$$
$$0 = [d'(x),R]_d,\eqno(6.23)$$
$$ a(x) = [d'(x), Q(x)]_d.\eqno(6.24)$$

 Relation (6.23) is automatically satisfied and, when $d'$ is diagonal,
(6.24) forces $a=0$.

 These algebraic constraints are 
essentially the ones of Theorems 1 and 2 of [dP1]. Indeed, their theorem 
1 gives an ansatz for a Lax pair, which contains our solution: they have three 
constraints, labelled there (3.7), (3.8), and (3.9); they then particularise  their
ansatz in Theorem 2 to what is in essence our case, with $d'$ diagonal;
 their conditions then particularise
to their (3.17), (3.18), (3.19). The first of their conditions follows automatically
from Weyl invariance; their second, (3.18), essentially tells us that $a=0$;
their third is condition (6.22). By choosing a suitable representation 
(which is strongly constrained by the conditions), they 
then ensure that these conditions can be satisfied.

The flow then has the Lax form (6.18) on the unreduced space $\tmgl$. Projecting 
to the reduced space, one  quotients out by the 
action of the automorphisms of the bundle, and so   omits the $d(x)$,
giving simply
$$\eqalign{  M^{-1}\dot M &= (\omega\phi^0_{Gl(N)})_d,\cr
(\dot \phi^0_{Gl(N)}) &= [(\omega \phi^0_{Gl(N)})^0_{od}
+ \omega \phi^0_{Gl(N)})^{\rm cst}_{od} ,(\phi^0_{Gl(N)})],\cr}
\eqno (6.25)$$
which is precisely the flow of the Hitchin system. In particular, one has a full set of 
commuting flows.

{\it The Bordner-Corrigan-Sasaki embedding}

  We keep our representation space $V= \bbc^N$ 
of a sum of weight 
spaces $\bbc_w$, invariant under the Weyl group, and still have our 
embedding $\hat\Xi$, turning our $H$-bundles $E(h)$ into bundles $E_{Gl(N)}=
\hat\Xi(E(h))
=\oplus_i L_{w_i(h)}$.
We now extend the embedding $\hat\Xi$ of the space of $H$ bundles to $\mred$
in a different way, corresponding to the Lax pairs of 
Bordner-Corrigan-Sasaki [BCS2]. This involves the construction of a different 
$Gl(N)$ Higgs field. 

 We first define some sections of $End(E_{Gl(N)})$
associated to sections of the bundles $L_{\alpha(h)}$. We note that 
to each root $\alpha$, we have an associated reflection $R_\alpha (v)
= v- <\hat\alpha, v> \alpha$ of the Lie algebra $\gh$, and in turn a permutation
of the weight spaces $\bbc_{w }$, which can be represented by a matrix
$(s_\alpha)_{w,w'}\in Gl(V)$, where, as usual we index the entries of the matrix 
by the weights
 themselves. Note that the non-zero entries of $s_\alpha$
must have    $w-w' = n\alpha$ for some integer $n$.
  For a section $f$ of $L_{\alpha(h)}$, we define 
$$\tilde s_\alpha(f)_{w,w'} =  \sum_n(s_\alpha)_{w,w'}\delta_{w-w',n \alpha} 
n\cdot I^*_nf.\eqno (6.26) $$

Represent the section $\phi$, which is
a section  of the associated bundle 
$E( (\oplus_\alpha\bbc_\alpha)\oplus\gh)$ 
by  $((\phi_\alpha), \phi_\gh)$.  We define the
 corresponding section $\phi_{Gl(N)}$   is given by  
$$\phi_{Gl(N)} =  ((\sum_\alpha \tilde s_\alpha(\phi_\alpha))
+\xi(\phi_\gh) ) .\eqno (6.27)$$

This section $\phi_{Gl(N)}$ has poles not only at the origin, but also,
for its $w, w+n\alpha$ components, at   the $n$-th roots of unity in the 
curve $\Sigma$. The moduli space of $Gl(N)$ Higgs pairs must be chosen accordingly.

The Hamiltonian is again a multiple of the Hamiltonian given
by pairing $tr(\phi_{Gl(N)}^2)$ with our standard cocycle $\omega$ of (5.7);
as above, the flow is given by (6.15)
$$\eqalign{  M^{-1}\dot M &= (\omega\phi^0_{Gl(N)})_d\cr
(\dot \phi^0_{Gl(N)})  &= a +[(\omega \phi^0_{Gl(N)})^0_{od},(\phi^0_{Gl(N)})]
+ 
[\omega \phi^0_{Gl(N)})^{\rm cst}_{od}, (\phi^0_{Gl(N)})]_{d} ,\cr
&= a +[(\omega \phi^0_{Gl(N)})^0_{od}
+ (\omega \phi^0_{Gl(N)})^{\rm cst}_{od},(\phi^0_{Gl(N)}) ]- 
[(\omega \phi^0_{Gl(N)})^{\rm cst}_{od}, (\phi^0_{Gl(N)}) ]_{od} ,\cr}
\eqno (6.28)$$
with again the equation(6.16) for $a$. 

Explicitly, one has 
$$(\phi^0_{Gl(N)})_{w,w'} = \sum_\alpha\sum_n(s_\alpha)_{w,w'}
\delta_{w-w',n \alpha} n \rho^0(\alpha(x), nz),$$
$$((\omega \phi^0_{Gl(N)})^0_{od}
+ (\omega \phi^0_{Gl(N)})^{\rm cst}_{od})_{w,w'} = \sum_\alpha\sum_n(s_\alpha)_{w,w'}
\delta_{w-w',n \alpha} {\partial\rho^0\over \partial x}(\alpha(x),nz).$$

One again wants to use the action of the diagonal subgroup to write (6.28)
 as a Lax pair. We take $d'(x)$, to be diagonal, and get
$$\eqalign{  M^{-1}\dot M &= (\omega\phi^0_{Gl(N)})_d\cr
(\dot \phi^0_{Gl(N)}) &= [(\omega \phi^0_{Gl(N)})^0_{od}
+ \omega \phi^0_{Gl(N)})^{\rm cst}_{od} + d'(x),(\phi^0_{Gl(N)})],\cr}.
\eqno (6.29)$$

In this case, the appropriate  diagonal terms are given in [BCS2]:
$$d'(x)_{w,w} = \sum_\alpha (s_\alpha)_{w,w }
  {\partial\rho^0\over \partial x}(\alpha(x),0).\eqno (6.30)$$

In this case, there is no constraint on the representation; indeed, Bordner,
Corrigan and Sasaki create a "universal Lax pair" within the algebra
$\bbc(H)\otimes \bbc[W]$ created by tensoring the group algebra of the Weyl group
with the function field of $H$; the product must be suitably defined, but
corresponds roughly to representing the Weyl group as acting by reflections
on a sum of weight spaces, and  the group $H$ as acting by diagonal matrices.
One then represents this algebra into a sum of weight spaces, and the Lax pair
gets embedded into $Gl(N)$. This  is what is given above. It would be interesting to do the geometry of
bundles directly within this algebra.

\bigskip
{\bf Bibliography}

\item{[BCS1]} A.J. Bordner, E. Corrigan and R. Sasaki, {\it Calogero-Moser
 models: I. A new formulation},  Progr. Theoret. Phys. 100 (1998), no. 6,
1107--1129. 
\item{[BCS2]} A.J. Bordner, E. Corrigan and R. Sasaki, {\it Generalized
 Calogero-Moser models and universal Lax pair operators}, hep-th/9905011
\item{[dP1]} E. d'Hoker and D.H. Phong, 
{\it Calogero-Moser Lax pairs with spectral parameter for general 
Lie algebras},  
Nuclear Phys. B {\bf 530} (1998), no. 3, 537--610 
\item{[dP2]} E. d'Hoker and D.H. Phong, 
{\it Spectral curves for
super-Yang-Mills with adjoint hypermultiplet for general simple Lie algebras}, 
Nuclear Phys. B {\bf 534} (1998), no. 3, 697--719
\item{[dP3]} E. d'Hoker and D.H. Phong, {\it Calogero-Moser and Toda
systems for twisted and untwisted affine Lie algebras},  
Nuclear Phys. B {\bf 530} (1998), no 3, 611-640
\item{[Do]} R. Donagi, {\it Seiberg-Witten integrable systems},  
Proc. Sympos. Pure Math., 62,
Part 2, Amer. Math. Soc., Providence, RI, 1997, 3--43
\item{[K]} I. M. Krichever, {\it Elliptic solutions of the 
Kadomtsev-Petviashvili equation
and integrable systems of particles}, Funct. Anal. Appl 14, 282-290 (1980)

\item{[Lo]} E. Looijenga,    {\it Root systems and elliptic curves}, Inv. 
Math. 38(1976),17-32 and {\it Invariant Theory for generalized root 
systems}, Inv. Math. 61,1-32(1980)

\item{[Ma]} E. Markman, {\it Spectral curves and integrable
systems}, Compositio Math. 93, 255-290  (1994).

\item{[OP]} M. A. Olshanetsky and A. M. Perelomov,
{\it Completely  integrable Hamiltonian systems connected 
with semisimple Lie algebras}, Inventiones Math 37 93-108 (1976)

\item{[OP2]} M. A. Olshanetsky and A. M. Perelomov, 
{\it Classical integrable finite-dimensional systems related 
to Lie algebras}, Phys. Rep. 71C, 313-400 (1981)

\ninerm 
\bigskip 
\bigskip 
 
\line{   J. C. Hurtubise\hfil E. Markman} 
\line{Centre de Recherches Math\'ematiques\hfil Department of Mathematics}
\line{Universit\'e de Montr\'eal \hfil   University of Massachusetts}
\line {and Department of Mathematics\hfil Amherst} 
\line{McGill University\hfil  
email: markman@math.umass.edu} 
\line{email: hurtubis@crm.umontreal.ca \hfil }

\end